\documentclass[12pt]{article}
\usepackage{amssymb,amsmath}
\usepackage{cases}
\usepackage{amsfonts}
\usepackage{cite,color,xcolor}
\usepackage[left=2.0cm,right=2.0cm,top=2.0cm,bottom=2.0cm]{geometry}
\usepackage[colorlinks,citecolor=blue,urlcolor=blue]{hyperref}

\newtheorem{theorem}{Theorem}[section]

\newtheorem{lemma}{Lemma}[section]

\newtheorem{remark}{Remark}[section]

\newcommand{\bal}{\begin{align}}
\newcommand{\bbal}{\begin{align*}}
\newcommand{\beq}{\begin{equation}}
\newcommand{\eeq}{\end{equation}}
\newcommand{\bca}{\begin{cases}}
\newcommand{\eca}{\end{cases}}

\newcommand{\pa}{\partial}
\newcommand{\fr}{\frac}

\newcommand{\de}{\delta}

\newcommand{\ep}{\varepsilon}
\newcommand{\dd}{\mathrm{d}}

\newcommand{\R}{\mathbb{R}}

\begin{document}
\title{Ill-posedness for the Burgers equation in Sobolev spaces}

\author{Jinlu Li$^{1}$, Yanghai Yu$^{2,}$\footnote{E-mail: lijinlu@gnnu.edu.cn; yuyanghai214@sina.com(Corresponding author); mathzwp2010@163.com} and Weipeng Zhu$^{3}$\\
\small $^1$ School of Mathematics and Computer Sciences, Gannan Normal University, Ganzhou 341000, China\\
\small $^2$ School of Mathematics and Statistics, Anhui Normal University, Wuhu 241002, China\\
\small $^3$ School of Mathematics and Big Data, Foshan University, Foshan, Guangdong 528000, China}

\date{\today}

\maketitle\noindent{\hrulefill}

{\bf Abstract:} In this paper, we considered the Cauchy problem for the Burgers equation and proved that the problem is ill-posed in Sobolev spaces $H^s$ with $s\in[1,\fr32)$.


\vskip0mm\noindent{\hrulefill}

\section{Introduction}\label{sec1}
\subsection{The Concept of Well-posedness}
A Cauchy problem
\begin{align*}
\pa_tf=F(f), \quad
f(0,x)=f_0(x)
\end{align*}
is said to be Hadamard well-posed in a Banach space $X$ if for any data $f_0\in X$ there exists $T>0$ and a unique solution in the space $\mathcal{C}([0,T),X)$ which depends continuously on the data. In particular, solutions describe continuous curves in $X$ at least for a short time. The problem is said to be ill-posed in X if it is not well-posed in the above sense. Based on the definition of well-posedness, there are at least three types of ill-posedness were studied in the literature: nonexistence, non-uniqueness, and discontinuous dependence on the data. In this paper we are interested in discontinuity with respect to the data.
\subsection{The Burgers equation}
The Burgers equation with fractional dissipation is written as
\begin{align}\label{fb}
\begin{cases}
\pa_tu+uu_x+\Lambda^{\gamma}u=0, &\quad (t,x)\in \R^+\times\R,\\
u(0,x)=u_0(x), &\quad x\in \R,
\end{cases}
\end{align}
where $\gamma\in[0,2]$ and the fractional power operator $\Lambda^{\gamma}$ is defined by Fourier multiplier with the symbol $|\xi|^{\gamma}$
\begin{eqnarray*}
  \Lambda^{\gamma}u(x)=\mathcal{F}^{-1}\big(|\xi|^{\gamma}\mathcal{F}u(\xi)\big).
\end{eqnarray*}
The Burgers equation \eqref{fb} with $\gamma=0$ and $\gamma=2$ has received an extensive amount of attention since the studies by Burgers in the 1940s. If $\gamma=0$, the equation is perhaps the most basic example of a PDE evolution leading to shocks. If $\gamma=2$, it provides an accessible model for studying the interaction between nonlinear and dissipative phenomena. Kiselev et al. \cite{Kiselev} gave a complete study for general $\gamma\in[0,2]$ for the periodic case. In particular, for the case $\gamma=1$, they proved the global well-posedness of the equation in the critical Hilbert space $H^{\fr12}(\mathbb{T})$ by using the method of modulus of continuity. Subsequently, Miao-Wu \cite{miao2009} proved the global well-posedness of the critical Burgers equation in critical Besov spaces $B^{1/p}_{p,1}(\R)$ with $p\in[1,\infty)$ with the help of Fourier localization technique and the method of modulus of continuity. For more results on the fractional Burgers equation and dispersive perturbations of Burgers equations, we refer the readers to see \cite{Alibaud,Dong,Karch,Linares,Molinet} and the references therein. We should mention that Molinet et al. \cite{Molinet} proved that the Cauchy problem for a class of dispersive perturbations of Burgers equations is locally well-posed in $H^s(\R)$.

In this paper, we focus on the well-posedness problem of the following Burgers equation.
\begin{align}\label{b}
\begin{cases}
\pa_tu+uu_x=0, &\quad (t,x)\in \R^+\times\R,\\
u(0,x)=u_0(x), &\quad x\in \R,
\end{cases}
\end{align}
Roughly speaking, \eqref{b} can be viewed as the simplest in the family of partial differential equations modeling the Euler and Navier-Stokes equation nonlinearity. The local well-posedness of the Burgers equation \eqref{b} for data in $H^s(\R)$ with any $s>3/2$ can be proved by combining the Sobolev embedding $H^{s-1}(\R)\hookrightarrow L^\infty(\R)$ and the classical energy estimate
$$\|u\|_{H^s}\leq \|u_0\|_{H^s}\exp\Big(C\int_0^t\|u_x\|_{L^\infty}\dd \tau\Big).$$
Also, in \cite{T} it is obtained that the solution map is continuous dependence while not uniformly continuous dependence on initial data for the Burgers equation \eqref{b} in the same space $H^s(\R)$ with $s>3/2$. For the endpoint case, Linares et al. in \cite{Linares} proved that the Cauchy problem for \eqref{b} is ill-posed in $H^{3/2}(\R)$, where the key point is that the available local well-posedness theory in $H^s(\R)$ with any $s>3/2$ have be used, for more details see Remark 1.6. Using the idea developed in \cite{Linares}, Guo et al. in \cite{Guo2019} to prove the ill-posedness for the Camassa-Holm equation in the critical Sobolev space $H^{3/2}(\R)$ and even in the Besov space $B^{1+1/p}_{p,r}(\R)$ with $r>1$. Precisely speaking, their main idea is to construct a blow-up smooth solution $u(t)\in \mathcal{C}([0,T^*),B^{1+1/p}_{p,r})\cap \mathcal{C}([0,T^*),H^2)$ such that
\bbal
\|u_0\|_{B^{1+\frac1p}_{p,r}}\leq \ep \quad\text{and} \quad \lim_{t\uparrow T^*}\|u(t)\|_{\mathrm{Lip}}=\infty,\quad T^*<\ep.
\end{align*}
then using the blow-up result and the following inequality
\bbal
\|u(t)\|_{\mathrm{Lip}}\leq C\|u(t)\|_{H^2}\leq C\big(1+\|u_0\|_{H^2}\big)\exp\exp\Big(C\int^t_0\|u_x\|_{B^0_{\infty,\infty}}\dd \tau\Big).
\end{align*}
they deduce that $\lim\limits_{t\uparrow T^*}\|u(t)\|_{B^1_{\infty,\infty}}=\infty$ which in turn implies $\lim\limits_{t\uparrow T^*}\|u(t)\|_{B^{1+\frac1p}_{p,r}}=\infty$. 

It remains nevertheless an interesting issue to prove that the Cauchy problem for the Burgers equation \eqref{b} in the Sobolev space $H^s(\R)$ with $s<\fr32$ is ill-posed. However, due to the absence of the embedding $H^s(\R)\hookrightarrow B^1_{\infty,\infty}(\R)$ for $s<3/2$, the method in \cite{Guo2019,Linares} is invalid for $s<3/2$.
In this paper, we shall develop a new method to study this problem and give a partial answer.
\subsection{Main Result}
Now let us state our main ill-posedness result of this paper.
\begin{theorem}\label{th}
Let $1\leq s<\frac32$. For any $\de>0$, there exists initial data satisfying
\bbal
\|u_0\|_{H^s}\leq \de,
\end{align*}
such that a solution $u(t)\in \mathcal{C}([0,T_0];H^s)$ of the Cauchy problem \eqref{b} satisfies
\bbal
\|u(T_0)\|_{H^s}\geq \frac1\de \quad \mathrm{for} \ \mathrm{some} \quad 0<T_0<\delta.
\end{align*}
\end{theorem}
\begin{remark}\label{re1}
Theorem \ref{th} indicates that the solutions of \eqref{b} with arbitrarily short time which initially have an arbitrarily small $H^s$-norm that grows arbitrarily large. This result shows the ill-posedness of \eqref{b} in $H^s(\R)$ with $1\leq s<\frac32$ in the sense that the solution map $u_0\in H^s\mapsto u\in H^s$ is discontinuous with respect to the initial data.
\end{remark}
{\bf Strategies to Proof.}\quad
We shall outline the main ideas in the proof of Theorem \ref{th}.
\begin{itemize}
  \item Firstly, we construct an explicit example for initial data $u_0$, where the norm $\|u_0\|_{H^s}$ is sufficiently small while $\|u'_0\|_{L^\infty}$ can be large enough.
  \item Secondly, we express the solution to the Burgers equation \eqref{b} by exploring fully the properties of the flow map and give the explicit blow-up time $T^*$.
  \item Lastly, we mainly observe that the transport term does cause growth of the $L^2$-norm of $u_x$ as $t$ tends to $T^*$. Precisely speaking, we estimate the $L^2$-norm of $u_x$ over $(-\psi(t,q_0),\psi(t,q_0))$ and obtain that its lower bound can be arbitrarily large as $t$ tends to $T^*$.
\end{itemize}
{\bf The structure of the paper.}\quad
In Section \ref{sec2} we provide several key Lemmas. In Section \ref{sec3} we present the proof of Theorem \ref{th}.

Let us complete this section with some notations we shall use throughout this paper.\\
{\bf Notations.}\quad The notation $A\leq a\wedge b$ means that $A \leq a$ and $A \leq b$. $a\approx b$ means $C^{-1}b\leq a\leq Cb$ for some positive harmless constants $C$. Given a Banach space $X$, we denote its norm by $\|\cdot\|_{X}$. For $I\subset\R$, we denote by $\mathcal{C}(I;X)$ the set of continuous functions on $I$ with values in $X$.
For all $f\in \mathcal{S}'$, the Fourier transform $\widehat{f}$ is defined by
$$
\widehat{f}(\xi)=\int_{\R}e^{-ix\xi}f(x)\dd x \quad\text{for any}\; \xi\in\R.
$$
For $s\in\R$, the nonhomogeneous Sobolev space $H^s(\R)$ is defined by its usual norm
\bbal
\|f\|^2_{H^s}=\int_{\R}(1+\xi^2)^s|\widehat{f}(\xi)|^2\dd \xi.
\end{align*}
\section{Preliminary}\label{sec2}
In the section, we make some preparations for the proof of the main theorem.
\subsection{Key Example for Initial Data}
Firstly, we construct an explicit example as follows. Set $$u_0(x):=p_0\big(e^{-|x+q_0|}-e^{-|x-q_0|}\big),$$
where two positive numbers $p_0$ and $q_0\in(0,1)$ will be fixed later.

It is easy to check that $u_0(x)$ is an odd function. Furthermore, we can deduce that the following result holds:
\begin{lemma}\label{le1} For every $q_0\in(0,1)$ and $s\in(\fr12,\fr32)$, there exists $C=C_s>0$ such that
\begin{align}\label{m}
C^{-1} p_0q_0^{3/2-s}\leq\|u_0\|_{H^s}\leq C p_0q_0^{3/2-s}.
\end{align}
\end{lemma}
{\bf Proof.}\quad The proof essentially follows that of Lemma 3.1 in \cite{Byers} or Proposition 1 in \cite{Himonas}. For the sake of readability, we sketch the proof here. Defining the function
\bbal
f(x):=e^{-|x+q_0|}-e^{-|x-q_0|},
\end{align*}
then using the fact $\widehat{e^{-|x|}}(\xi)=2 /\left(1+\xi^{2}\right)$, we have
\bbal
\widehat{f}(\xi)=\frac{2(e^{iq_0\xi}-e^{-iq_0\xi})}{1+\xi^2}=\frac{4i\sin(q_0\xi)}{1+\xi^2}.
\end{align*}
Using the definition of the $H^s$-norm and the change of variable setup $y=q_0\xi$, we have
\begin{align*}
\left\|f\right\|_{H^{s}(\mathbb{R})}^{2}&=16\int_{\mathbb{R}}\left(1+\xi^{2}\right)^{s-2} \sin ^{2}(q_0\xi) \dd \xi\nonumber\\
&\geq32\Big(1+\fr{\pi^2}{q_0^2}\Big)^{s-2}\int_{0}^{\pi/q_0}\sin ^{2}(q_0\xi) \dd \xi\nonumber\\
&=32\Big(1+\fr{\pi^2}{q_0^2}\Big)^{s-2}\cdot\fr1{q_0}\int_{0}^{\pi}\sin ^{2}y \dd y\nonumber\\
&\geq 16\pi\left(1+\pi^{2}\right)^{-3/2} q_0^{3-2s},
\end{align*}
where we have used $s\in\left(\frac{1}{2}, \frac{3}{2}\right)$ and $q_0 \in(0,1)$. This proves the lower bound.

To get the upper bound, we split the domain of integration as
\begin{align*}
\left\|f\right\|_{H^{s}(\mathbb{R})}^{2}&=32\int_{0}^{\infty}\left(1+\xi^{2}\right)^{s-2} \sin ^{2}(q_0\xi) \dd \xi\nonumber\\
&=32\Big(\int_{0}^{1/q_0}+\int_{1/q_0}^{\infty}\Big)\left(1+\xi^{2}\right)^{s-2} \sin ^{2}(q_0\xi) \dd \xi\nonumber\\
&=:32(I_1+I_2).
\end{align*}
Due to the simple fact $\sin^{2}(q_0\xi)\leq |q_0\xi|^2\wedge1$, we have
\begin{align*}
&I_1\leq32q_0^2\int_{0}^{1/q_0}\xi^{2s-2} \dd \xi\leq \Big(\frac{32}{2s-1}\Big)q_0^{3-2s},\nonumber\\
&I_2\leq32\int_{1/q_0}^{\infty}\xi^{2s-4} \dd \xi\leq \Big(\frac{32}{3-2s}\Big)q_0^{3-2s},
\end{align*}
which completes the proof of Lemma \ref{le1}.
\subsection{Existence and Blow-up criterion}
\begin{lemma}\label{le2} For every $s\in[1,\fr32)$, there exists a solution $u\in \mathcal{C}([0,T^*);H^s)\cap L^\infty([0,T^*);\mathrm{Lip})$ for the Burgers equation \eqref{b}, where $T^*<\infty$ is the maximal time for initial data $u_0$.

Furthermore, we have
\bbal
\lim_{t\uparrow T^*}\big(\|u(t)\|_{H^s}+\|u(t)\|_{\mathrm{Lip}}\big)=+\infty\quad\Leftrightarrow\quad\lim_{t\uparrow T^*}\|u(t)\|_{\mathrm{Lip}}
=+\infty.
\end{align*}
\end{lemma}
{\bf Proof.}\quad Easy computations give that
\bal\label{u}
u'_0(x)=
\begin{cases}
-p_0(e^{-q_0}-e^{q_0})e^{x}, &\;\text{if}\; x\in(-\infty,-q_0),\\
-p_0e^{-q_0}(e^x+e^{-x}), &\;\text{if}\; x\in(-q_0,q_0),\\
-p_0(e^{-q_0}-e^{q_0})e^{-x}, &\;\text{if}\; x\in(q_0,+\infty),
\end{cases}
\end{align}
from which and Lemma \ref{le1}, we get that $u_0\in H^s\cap \mathrm{Lip}$.

Following the proof of Lemma 2.4 in \cite{d1}, we can get
\begin{align*}
&\|u(t)\|_{H^s} \leq\left\|u_{0}\right\|_{H^s} \exp\Big(C \int_{0}^{t}\|u(\tau)\|_{\mathrm{Lip}} \dd \tau\Big)
\end{align*}
and
\begin{align*}
&\|u(t)\|_{\mathrm{Lip}} \leq\left\|u_{0}\right\|_{\operatorname{Lip}} \exp\Big(C \int_{0}^{t}\|u(\tau)\|_{\mathrm{Lip}} \dd \tau\Big).
\end{align*}
This is enough to complete the proof of Lemma \ref{le2}.
\subsection{The Equation Along the Flow}
Given a Lipschitz velocity field $u$, we may solve the following ODE to find the flow induced by $u$:
\begin{align}\label{ode}
\quad\begin{cases}
\frac{\dd}{\dd t}\psi(t,x)=u(t,\psi(t,x)),\\
\psi(0,x)=x,
\end{cases}
\end{align}
which is equivalent to the integral form
\bbal
\psi(t,x)=x+\int^t_0u(\tau,\psi(\tau,x))\dd \tau.
\end{align*}
Furthermore, we get from \eqref{b} that
\bbal
\frac{\dd}{\dd t}u(t,\psi(t,x))&=u_{t}(t,\psi(t,x))+u_{x}(t,\psi(t,x))\frac{\dd}{\dd t}\psi(t,x)=0,
\end{align*}
which means that
\bal\label{l1}
u(t,\psi(t,x))=u_0(x),\quad \text{namely},\quad u(t,x)=u_0(\psi^{-1}(t,x)).
\end{align}
Thus we can give the explicit expression of the flow as:
\bal\label{l2}
\psi(t,x)=x+tu_0(x).
\end{align}
Let $y=\psi(t,x)$, then we have
\bal\label{z2}
\psi^{-1}(t,y)=y-tu_0(\psi^{-1}(t,y))=y-tu(t,y).
\end{align}
Differentiating \eqref{b} with respect to space variable $x$, we find
\bbal
u_{tx}+uu_{xx}+(u_x)^2=0.
\end{align*}
Combining the above and \eqref{ode}, we obtain
\bbal
\frac{\dd}{\dd t}u_x(t,\psi(t,x))&=u_{tx}(t,\psi(t,x))+u_{xx}(t,\psi(t,x))\frac{\dd}{\dd t}\psi(t,x),\\
&=u_{tx}(t,\psi(t,x))+u_{xx}(t,\psi(t,x))u(t,\psi(t,x))\\
&=-(u_x)^2(t,\psi(t,x)),
\end{align*}
which reduces to
\bal\label{l3}
u_x(t,\psi(t,x))=\frac{1}{t+\frac{1}{u'_0(x)}}.
\end{align}
We should mention that the above can also be deduced from \eqref{l1} and \eqref{l2}.

According to the definition of $u_0$, we can deduce that
$$T^*=-\frac{1}{u'_0(q^-_0)}\in \Big(\frac{1}{2p_0},\frac{1}{p_0}\Big).$$
Because the velocity field is Lipschitz, then we get that for $t\in[0,T^*)$
\bbal
\psi_x(t,x)=\exp\Big(\int^t_0u_x(\tau,\psi(\tau,x))\dd \tau\Big)>0.
\end{align*}
This shows that $\psi(t,\cdot)$ is an increasing diffeomorphism over $\R$, that is,  for all $x,y\in \R,$ there holds that $\psi(t,x)<\psi(t,y)$ if $x< y$.

\section{Proof of Main Theorem}\label{sec3}
In this section, we prove Theorem \ref{th}.

{\bf Proof of Theorem \ref{th}}\;  By the definition of $u'_0(x)$ and \eqref{l3}, we know that $u_x(t,\psi(t,x))$ is continuous in $[0,T^*)\times (-q_0,q_0)$. We should emphasize that $u_x(t,x)$ is discontinuous in $[0,T^*)\times \R$, but we can claim that $u_x(t,x)$ is continuous in $[0,T^*)\times (-\psi(t,q_0),\psi(t,q_0))$. In fact, for any $x,y\in (-\psi(t,q_0),\psi(t,q_0))$, we have $\psi^{-1}(t,x),\psi^{-1}(t,y)\in(-q_0,q_0)$.  Moreover, we deduce from \eqref{z2} that  for $t\in[0,T^*)$
\bal\label{l4}
|\psi^{-1}(t,x)-\psi^{-1}(t,y)|&\leq |\pa_x\psi^{-1}(t,x)||x-y|\nonumber\\
&\leq|x-y| \big(1+T^*\|u_x\|_{L^\infty_t(L^\infty)} \big).
\end{align}
Also, it follows from \eqref{z2} and \eqref{b} that
\bal\label{z4}
|\psi^{-1}(t,x)-\psi^{-1}(s,x)|&=|tu(t,x)-su(s,x)| \nonumber
\\&\leq \|u_0\|_{L^\infty}|t-s|+T^*\Big|\int^t_s\|\pa_tu(\tau,\cdot)\|_{L^\infty}\dd \tau\Big| \nonumber
\\&\leq C\big(1+T^*\|u_x\|_{L^\infty_t(L^\infty)} \big)|t-s|.
\end{align}
Thus, we obtain from \eqref{u} and \eqref{l4}-\eqref{z4} for $s,t\in[0,T^*)$ and $x,y\in (-\psi(t,q_0),\psi(t,q_0))$
\bbal
|u_x(t,x)-u_x(s,y)|&\leq|u_x(t,x)-u_x(t,y)|+|u_x(t,y)-u_x(s,y)|\\
&\leq|u_x(t,\psi(t,\psi^{-1}(t,x)))-u_x(t,\psi(t,\psi^{-1}(t,y)))|\\
&~~~~+|u_x(t,\psi(t,\psi^{-1}(t,y)))-u_x(s,\psi(s,\psi^{-1}(s,y)))|\\
&\to0\quad\text{as}\quad (t,x)\to (s,y).
\end{align*}
By the Burgers equation $\pa_tu=-uu_x$, we can deduce that $\pa_tu(t,x)$ is continuous in $[0,T^*)\times (-\psi(t,q_0),\psi(t,q_0))$. That is $u(t,x)\in \mathcal{C}^1([0,T^*)\times(-\psi(t,q_0),\psi(t,q_0)))$. Furthermore, one has
\bbal
u_{xx}(t,\psi(t,x))=\frac{u''_0(x)}{\big(1+tu'_0(x)\big)^3}.
\end{align*}
The similar argument shows that $u_x(t,x)\in \mathcal{C}^1([0,T^*)\times(-\psi(t,q_0),\psi(t,q_0)))$.

For notational convenience we now set
\bbal
&\widetilde{m}(t):=u_x(t,\psi(t,0))=\frac{1}{t+\frac{1}{u'_0(0)}}.
\end{align*}
Therefor, we have
\bal\label{l}
u_x(t,\psi(t,x))\leq \widetilde{m}(t)\quad\text{for all}\; x\in(-q_0,q_0).
\end{align}
Set $w(t,x):=u_x(t,x)$, then we obtain from \eqref{b}
\bbal
\pa_tw+\pa_x(uw)=0,
\end{align*}
which implies that for $(t,x)\in[0,T^*)\times(-\psi(t,q_0),\psi(t,q_0))$
\bal\label{y}
\pa_t(w^2)+\pa_x(uw^2)+\pa_xuw^2=0.
\end{align}
Integrating \eqref{y} with respect to space variable $x$ over $[-\psi(t,q_0),\psi(t,q_0)]$, we have
\bal\label{y1}
\int_{|x|\leq \psi(t,q_0)}\pa_t(w^2)\dd x+\int_{|x|\leq \psi(t,q_0)}\pa_x(uw^2)\dd x+\int_{|x|\leq \psi(t,q_0)}\pa_xuw^2\dd x=0.
\end{align}
As $u_0(x)$ is odd, the solution of Burgers equation satisfies $u(t,x)=-u(t,-x)$, which tells us that $w(t,x)=w(t,-x)$. Thus we have
\bal\label{y2}
\int_{|x|\leq \psi(t,q_0)}\pa_t(w^2)\dd x=\frac{\dd}{\dd t}\int_{|x|\leq \psi(t,q_0)}w^2\dd x-2u(t,\psi(t,q_0))w^2(t,\psi^-(t,q_0)),
\end{align}
and
\bal\label{y3}
\int_{|x|\leq \psi(t,q_0)}\pa_x(uw^2)\dd x=2u(t,\psi(t,q_0))w^2(t,\psi^-(t,q_0)).
\end{align}
Inserting \eqref{y2} and \eqref{y3} into \eqref{y1} yields
\bal\label{y4}
\frac{\dd}{\dd t}\int_{|x|\leq \psi(t,q_0)}w^2\dd x+\int_{|x|\leq \psi(t,q_0)}\pa_xu(t,x)w^2\dd x=0.
\end{align}
To simplify notation let $$A(t):=\int_{|x|\leq \psi(t,q_0)}w^2(t,x)\dd x\quad\text{for}\quad t\in[0,T^*),$$
combining \eqref{l}, then \eqref{y4} reduces to
\bbal
A'(t)=\int_{|x|\leq \psi(t,q_0)}-u_x(t,x)w^2\dd x\geq -\widetilde{m}(t)A(t).
\end{align*}
Solving the above differential inequality gives us that
\bbal
A(t)\geq A_0\exp\Big(\int^t_0-\widetilde{m}(\tau)\dd \tau\Big)= A_0\cdot \frac{\widetilde{m}(t)}{\widetilde{m}(0)}.
\end{align*}
which implies
\bal\label{y5}
\|w\|_{L^2}\geq A^{\frac12}_0\cdot \sqrt{\frac{\widetilde{m}(t)}{u'_0(0)}}.
\end{align}
Notice that
\bbal
A_0=\int_{|x|\leq q_0}\big(u'_0(x)\big)^2\dd x\approx  p^2_0q_0,
\end{align*}
and
\bbal
\lim_{t\uparrow T^*}\widetilde{m}(t)&=\frac{1}{-\frac{1}{u'_0(q_0^-)}+\frac{1}{u'_0(0)}}
=\frac{u'_0(q_0^-)u'_0(0)}{u'_0(q_0^-)-u'_0(0)},
\end{align*}
combining the above and \eqref{y5} yields\bbal
\lim_{t\uparrow T^*}\|u(t)\|_{H^1}&\geq \lim_{t\uparrow T^*}\|w\|_{L^2}\\
&\geq Cp_0\sqrt{q_0}\sqrt{\frac{u'_0(q_0^-)}{u'_0(q_0^-)-u'_0(0)}}\\
&\geq C\frac{p_0\sqrt{q_0}}{1-e^{-q_0}}\\
&\approx Cp_0q_0^{-\fr12},
\end{align*}
where we have used that
$$u'_0(0)=-2p_0e^{-q_0}\quad\text{and}\quad u'_0(q_0^-)=-p_0(e^{-2q_0}+1)$$
and in the last step used
$$\frac{q_0}{2}\leq1-e^{-q_0}\leq q_0\quad\text{for}\quad q_0\in(0,1).$$
By Lemma \ref{le1}, one has
\bbal
\|u_0\|_{H^s}\leq c_1p_0q^{\frac32-s}_0\leq \delta \quad\text{and}\quad T^*\leq \frac{1}{p_0}\leq \delta,
\end{align*}
but
\bbal
\lim_{t\uparrow T^*}\|u(t)\|_{H^1}
&\geq c_2p_0q_0^{-\fr12}\geq \frac{1}{\delta^2},
\end{align*}
if some large $p_0$ and small $q_0$ is chosen. In fact, we can take $p_0$ sufficiently large to make $p_0\geq \frac1\delta$ and $q_0$ sufficiently small such that
$q_0\leq \big(\delta/(c_1p_0)\big)^{2/(3-2s)}\wedge c_2^2p_0^2\delta^4$.

Hence, we can choose $T_0\in[0,T^*)$ such that
\bbal
\|u(T_0)\|_{H^s}\geq C\|u(T_0)\|_{H^1}\geq \frac1\delta.
\end{align*}
This completes the proof of Theorem \ref{th}.

\section*{Acknowledgments} J. Li is supported by the National Natural Science Foundation of China (Grant No.11801090). Y. Yu is supported by the Natural Science Foundation of Anhui Province (No.1908085QA05). W. Zhu is partially supported by the National Natural Science Foundation of China (Grant No.11901092) and Natural Science Foundation of Guangdong Province (No.2017A030310634).

\end{document}